\documentclass{gtart_h}

\def\ifplaintex{\expandafter\ifx\csname documentclass\endcsname\relax}

\def\gtp{{\mathsurround=0pt\it $\cal G\mskip-2mu$eometry \&\ 
$\cal T\!\!$opology $\cal P\!$ublications}}  

\def\recd{{\small Received:\qua\receiveddate\ifx\reviseddate\relax
\else\qquad Revised:\qua\reviseddate\fi\par}} 

\def\lognumber#1{\def\thelognumber{#1}}
\def\volumenumber#1{\def\thevolumenumber{#1}}
\def\volumeyear#1{\def\thevolumeyear{#1}}
\def\papernumber#1{\def\thepapernumber{#1}}
\def\pagenumbers#1#2{\def\startpage{#1}\def\finishpage{#2}}
\def\published#1{\def\publishdate{#1}}

\def\received#1{\def\receiveddate{#1}}
\def\revised#1{\def\reviseddate{#1}}
\def\accepted#1{\def\accepteddate{#1}}

\long\def\asciiabstract#1{\long\def\theasciiabstract{#1}}

\let\\\par\let\thelognumber\relax\let\thevolumenumber\relax
\let\thepapernumber\relax\let\thevolumeyear\relax\let\startpage\relax
\let\finishpage\relax\let\publishdate\relax\let\receiveddate\relax
\let\reviseddate\relax\let\accepteddate\relax\let\theasciititle\relax
\let\theasciiauthors\relax
\let\theasciiabstract\relax

\let\theasciiemail\relax

\ifplaintex
\font\logobig=cmssbx10 scaled 3836
\font\logomed=cmssbx10 scaled 2557
\else
\font\logobig=cmssbx10 scaled 4200
\font\logomed=cmssbx10 scaled 2800
\fi

\long\def\makeagttitle{   
\count0=\startpage
\agt\hfill      
\hbox to 45truept{\vbox to 0pt{\vglue -13truept{\logomed A\kern -.37em{\logobig 
T}\kern -.38em G}\vss}\hss}
\break
{\small Volume \thevolumenumber\ (\thevolumeyear)
\startpage--\finishpage\nl
Published: \publishdate}

\vglue .25truein

{\parskip=0pt\leftskip 0pt plus
1fil\def\\{\par\smallskip}{\Large\bf\thetitle}\par\medskip} \vglue
0.05truein

{\parskip=0pt\leftskip 0pt plus 1fil\def\\{\par}{\sc\theauthors}
\par\medskip}%
 
\vglue 0.03truein 

{\small\leftskip 25truept\rightskip 25truept{\bf Abstract}\stdspace\theabstract

{\bf AMS Classification}\stdspace\theprimaryclass
\ifx\thesecondaryclass\relax\else; \thesecondaryclass\fi\par
{\bf Keywords}\stdspace \thekeywords\par}\vglue 7truept

}   

\ifplaintex
\font\phead=cmsl9 scaled 950
\font\pnum=cmbx10 scaled 913
\font\pfoot=cmsl9 scaled 950
\headline{\vbox to 0pt{\vskip -4.5mm\line{\small\phead\ifnum
\count0=\startpage ISSN 1472-2739 (on-line) 1472-2747 (printed)
\hfill {\pnum\folio}\else\ifodd\count0\def\\{ }%
\ifx\theshorttitle\relax\thetitle\else\theshorttitle\fi\hfill{\pnum\folio}
\else\def\\{ and }{\pnum\folio}\hfill\ifx\theshortauthors\relax\theauthors
\else\theshortauthors\fi\fi\fi}\vss}}
\footline{\vbox to 0pt{\vglue 0mm\line{\small\pfoot\ifnum\count0=\startpage
\copyright\ \gtp\hfill\else
\agt, Volume \thevolumenumber\ (\thevolumeyear)\hfill\fi}\vss}}
\else
\font=cmsl9 scaled 1050
\font\lnum=cmbx10 
\font=cmsl9 scaled 1050
\makeatletter
\def\@oddhead{{\small\lhead\ifnum\count0=\startpage ISSN 1472-2739 
(on-line) 1472-2747 (printed)\hfill {\lnum\number\count0}\else\ifodd\count0
\def\\{ }\ifx\theshorttitle\relax \thetitle \else\theshorttitle\fi\hfill
{\lnum\number\count0}\else\def\\{ and }{\lnum\number\count0}
\hfill\ifx\theshortauthors\relax 
\theauthors\else\theshortauthors\fi\fi\fi}}\def\@evenhead{\@oddhead}
\def\@oddfoot{\small\lfoot\ifnum\count0=\startpage\copyright\ \gtp\hfill\else
\agt, Volume \thevolumenumber\ (\thevolumeyear)\hfill\fi}
\def\@evenfoot{\@oddfoot}
\makeatother
\fi
\let\maketitlepage\makeagttitle
\let\makeshorttitle\maketitlepage
\let\maketitle\maketitlepage

\newwrite\gtoutfile
\long\gdef\makeheadfile{  
{\def\\{, }\def\s{ }
\immediate\openout\gtoutfile head.xxx
\immediate\write\gtoutfile{Proxy-for: \ifx\theasciiauthors\relax
\theauthors\else\theasciiauthors\fi\s<\ifx\theasciiemail\relax\theemail\else\theasciiemail\fi>}
\immediate\write\gtoutfile{\noexpand\\}
\immediate\write\gtoutfile{Authors: \ifx\theasciiauthors\relax
\theauthors\else\theasciiauthors\fi}
{\def\\{ }\immediate\write\gtoutfile{Title: \ifx\theasciititle\relax
\thetitle\else\theasciititle\fi}}
\immediate\write\gtoutfile{Subj-class: GT or SG, GR etc}
\immediate\write\gtoutfile{MSC-class: \theprimaryclass\ifx\thesecondaryclass\relax\else, \thesecondaryclass\fi}
\immediate\write\gtoutfile{Journal-ref: Algebr. Geom. Topol. \thevolumenumber\s
(\thevolumeyear) \startpage-\finishpage}
\immediate\write\gtoutfile{Comments: Published by Algebraic and
Geometric Topology at}
\immediate\write\gtoutfile{\s\s\s  http://www.maths.warwick.ac.uk/agt/AGTVol\thevolumenumber/agt-\thevolumenumber-\thepapernumber.abs.html}
\immediate\write\gtoutfile{\noexpand\\}
\immediate\write\gtoutfile{}
\ifx\theasciiabstract\relax
\immediate\write\gtoutfile{\theabstract}\else
\immediate\write\gtoutfile{\theasciiabstract}\fi
\immediate\write\gtoutfile{}
\immediate\write\gtoutfile{\noexpand\\}
\immediate\write\gtoutfile{}
\immediate\closeout\gtoutfile}}  

\def\maketitlepage{\makeagttitle\makeheadfile}
\let\makeshorttitle\maketitlepage
\let\maketitle\maketitlepage

\lognumber{48}
\volumenumber{5}
\volumeyear{2005}
\papernumber{48}
\pagenumbers{1197}{1222}
\received{9 August 2004} 
\revised{23 July 2005}
\accepted{14 March 2005}
\published{20 September 2005}

\usepackage{amsmath,amssymb,graphicx,psfrag}

\def\psfraga <#1,#2> #3#4{%
\psfrag {#3}{\smash{\rlap{\kern #1 \raise #2\hbox{#4}}}}}

\def\fref#1{\hyperlink{#1anchor}{\ref*{#1}}}
\def\figref#1{\hyperlink{#1anchor}{Figure~\ref*{#1}}}
\def\anchor#1{\noindent\hypertarget{#1anchor}{\smash{$\phantom{99}$}}}

\newcommand{\HF}{HF}

\newtheorem{theorem}{Theorem}[section]
\newtheorem{prop}[theorem]{Proposition}
\newtheorem{cor}[theorem]{Corollary}
\newtheorem{lemma}[theorem]{Lemma}

\theoremstyle{definition}

\newtheorem{defn}[theorem]{Definition}

\newcommand{\R}{\mathbb{R}}
\newcommand{\T}{\mathbb{T}}

\newcommand{\Z}{\mathbb{Z}}

\newcommand{\cm}{\cdot}

\newcommand\relspinc{\underline{\spinc}}

\newcommand\Filt{\mathcal F}

\newcommand\x{\mathbf x}
\newcommand\w{\mathbf w}
\newcommand\z{\mathbf z}
\newcommand\p{\mathbf p}
\newcommand\q{\mathbf q}
\newcommand\y{\mathbf y}

\newcommand\ModSphere{\ModFlow\left({\mathbb S}\longrightarrow 
\Sym^{g-1}(\Sigma_{1})\times \Sym^2(\Sigma_{2})\right)}
\newcommand\ModSpheres\ModSphere

\newcommand\CFa{\widehat{CF}}

\newcommand\HFp{\HFb}

\newcommand\HFm{\HF^-}

\newcommand\HFinf{HF^\infty}

\newcommand\HFa{\widehat{HF}}
\newcommand\HFb{HF^+}
\newcommand\gr{\mathrm{gr}}
\newcommand\Mas{\mu}
\newcommand\UnparModSp{\widehat \ModSp}
\newcommand\UnparModFlow\UnparModSp
\newcommand\Mod\ModSp

\newcommand{\spinc}{\mathfrak s}

\newcommand\ModMaps{\mathcal M}
\newcommand\ModSp\ModMaps

\newcommand\Ta{{\mathbb T}_{\alpha}}
\newcommand\Tb{{\mathbb T}_{\beta}}

\newcommand\spincrel\relspinc

\newcommand\CFK{CFK}
\newcommand\HFK{HFK}

\newcommand\CFKa{\widehat\CFK}

\newcommand\HFKa{\widehat\HFK}

\newcommand\BasePt{w}
\newcommand\FiltPt{z}

\newcommand\Dual{\mathcal D}
\newcommand\Duality\Dual

\begin{document}

\title{On knot Floer homology and cabling}
\author{Matthew Hedden}
\address{Department of
Mathematics, Princeton University\\Princeton, NJ 08544-1000, USA}
  
\email{mhedden@math.princeton.edu}

\begin{abstract}
This paper is devoted to the study of the knot Floer
homology groups $\HFKa(S^3,K_{2,n})$, where $K_{2,n}$ denotes the
$(2,n)$ cable of an arbitrary knot, $K$.  It is shown that for
sufficiently large $|n|$, the Floer homology of the cabled knot
depends only on the filtered chain homotopy type of $\CFKa(K)$.  A precise formula for this relationship is presented.   
In fact, the homology groups in the top $2$ filtration dimensions for
the cabled knot are
isomorphic to the original knot's Floer
homology group in the top filtration dimension. The results are
extended to $(p,pn\pm1)$ cables.  As an example we compute
 $\HFKa((T_{2,2m+1})_{2,2n+1})$ for all sufficiently large $|n|$, where $T_{2,2m+1}$ denotes the
$(2,2m+1)$-torus knot.
\end{abstract}

\asciiabstract{This paper is devoted to the study of the knot Floer
homology groups HFK(S^3,K_{2,n}), where K_{2,n} denotes the
(2,n) cable of an arbitrary knot, K.  It is shown that for
sufficiently large |n|, the Floer homology of the cabled knot
depends only on the filtered chain homotopy type of CFK(K).  A
precise formula for this relationship is presented.  In fact, the
homology groups in the top 2 filtration dimensions for the cabled
knot are isomorphic to the original knot's Floer homology group in the
top filtration dimension. The results are extended to (p,pn+-1)
cables.  As an example we compute HFK((T_{2,2m+1})_{2,2n+1}) for
all sufficiently large |n|, where T_{2,2m+1} denotes the
(2,2m+1)-torus knot.}

\maketitle

\section{Introduction}

In \cite{HolDisk}, Ozsv{\'a}th and Szab{\'o} introduced a collection of
abelian groups associated to closed oriented three-manifolds: given a
three-manifold $Y$ and Spin$^c$ structure $\spinc$, there are various Heegaard Floer homology
groups of $Y$:  $\HFa(Y,\spinc)$,
$\HFinf(Y,\spinc)$, $\HFp(Y,\spinc)$, and $\HFm(Y,\spinc)$.
In \cite{Knots} they subsequently showed that a knot $K\subset Y$
induces a filtration on the
chain complexes which compute these groups, see also \cite{RasmussenThesis}. In particular, the filtered chain homotopy
types of the filtered chain complexes were shown
to be topological invariants of the knot and the Spin$^c$
structure. This paper will deal  with the case $Y=S^3$ and
primarily with the simplest objects defined in \cite{Knots}, $\HFKa(K,i)$.
The notation here, as in the rest of this paper, agrees
whenever possible with that of 
\cite{Knots},\cite{AltKnots},\cite{GBmutant}, so that the lack of
reference to the three-manifold in $\HFKa(K,i)$ implies $Y=S^3$ and
the index $i$ refers to the level of the filtration induced on
$\HFa(S^3)$ by $K$. 
Specific definitions and relevant notation will be discussed in Section \ref{sec:Heegaards}.

\begin{figure}[ht!]\small\anchor{fig:projection}
\psfrag{+1}{$+1$}
\begin{center}
\includegraphics{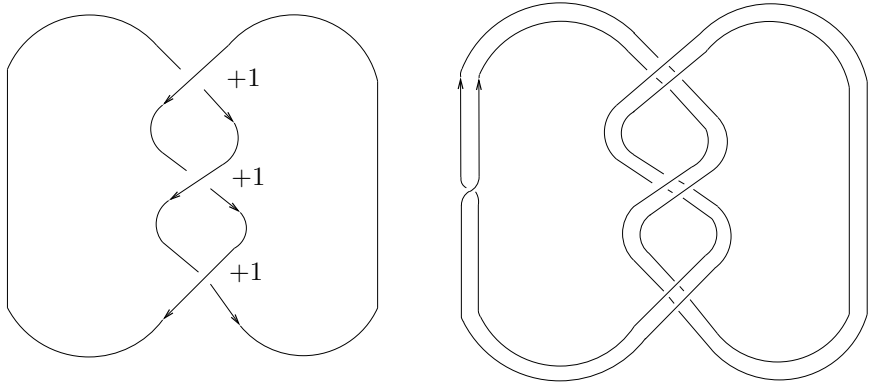}
\caption{\label{fig:projection} On the left is the right handed trefoil, with writhe +3.  On
  the right is the (2,7) cable of the trefoil.}
\end{center}
\end{figure}

Recall that the $(p,q)$ cable of a knot $K$, denoted $K_{p,q}$, is defined to be the
topological type of a knot supported on the boundary of a
tubular neighborhood of $K$ which is
linear with slope $p/q$ with
respect to the standard framing of this torus.  In other words, it is
a satellite knot which
winds $p$ times around the meridian of $K$ as it winds $q$ times
around a specified longitude.  This longitude is determined by the Seifert
framing for $K$. The knot which is 
cabled is sometimes called the companion knot (see \cite{Lickorish}
for more details).
 
Cabling a knot increases its complexity in some sense.
If one draws a projection for a knot and its $(p,pn+1)$
cable (where $n$ is the writhe of the original knot's projection), the number of crossings in the
latter projection will be $p^2$ times the number of crossings of the
original diagram plus $(p-1)$ (see \figref{fig:projection}).

The Alexander polynomials of a knot and its cables are related by the
following classical formula (which is a special case of a similar
formula holding for all satellites):
\begin{equation}
\label{eq:satpoly}
\Delta_{K_{p,q}}(t)=\Delta_{T_{p,q}}(t)\cm\Delta_K(t^p),
\end{equation}
where $T_{p,q}$ denotes the $(p,q)$ torus knot, and $\Delta_K(t)$ the
symmetrized Alexander polynomial of $K$ \cite{Lickorish}.
It is proved in \cite{Knots} that the following relationship holds
between the Euler characteristics of $\HFKa$ and the symmetrized
Alexander polynomial:
\begin{equation}
\label{eq:Euler}
\sum_{i} \chi \left(\HFKa(K,i)\right) \cm T^i =\Delta_K(T).
\end{equation}
It is therefore a natural question to ask how
Equation \eqref{eq:satpoly} manifests itself within $\HFKa$, and, more 
generally, how $\HFKa$ of a knot and its satellites are related.

We demonstrate some results in this direction. Before
stating the first theorem, recall from \cite{GBmutant} that
deg$\HFKa(K)$ denotes the largest integer $d>0$ for which
$\HFKa(K,d)\ne0$ and that \cite{genusbounds} identifies this
invariant with the Seifert genus of $K$.   Note that also
deg$\HFKa(T_{p,q})=\frac{(p-1)(q-1)}{2}$.     Let us denote the
filtration of $\CFa(S^3)$ induced by $K$ by $\Filt(K,j)$, so that we
have the sequence of inclusions:
$$ 0=\Filt(K,-i) \subseteq \Filt(K,-i+1)\subseteq \ldots \subseteq
\Filt(K,n)=\CFa(S^3),$$
\noindent with $\frac{\Filt(K,j)}{\Filt(K,j-1)}=\CFKa(K,j)$. The
following theorem will be proved in Section~\ref{sec:largen}.

\begin{theorem}
\label{thm:2ncable}
Let $K$ be a knot in $S^3$, and suppose $\mathrm{deg}~\HFKa(K)=d$. Then $\exists \
N>0$ such that $\forall \ n>N$ the following holds:
$$ {\rm deg}~\HFKa(K_{2,2n+1})=2d+n. $$
Furthermore, $\forall \ i \ge 0$ we have
$$ \HFKa_*(K_{2,2n+1},i)\cong 
\left\{\begin{array}{ll}

 H_{*+2(k-d)}(\Filt(K,k-d)) & {\text{for
 $i=2d+n-2k$}} \\

H_{*+2(k-d)+1}(\Filt(K,k-d)) & {\text{for
 $i=2d+n-2k-1$.}} \\

\end{array}
\right. $$
\end{theorem}

By the symmetry of $\HFKa$ under the involution on Spin$^c$
structures (Equation \eqref{eq:jsym} in Section \ref{sec:Heegaards}) the above result completely
determines $\HFKa(K_{2,2n+1})$.
Note that the information required above is more
than simply $\HFKa(K,i)$ for all $i$. $\HFKa(K)$ is the homology of an associated
graded of a filtered chain complex -- one needs to know
$H_*(\Filt(K))$ to fully exploit the theorem.   
Despite this additional requirement, the theorem is still a
powerful calculational tool.  For instance, in \cite{Goda} it is
shown that the Floer homology of (1,1) knots is combinatorial.  They show this by exhibiting a genus one Heegaard
diagram for a (1,1) knot.  Since the differentials in these cases can
be computed combinatorially via the Riemann mapping theorem, $\HFKa$
of $(2,2n+1)$ cables will be
given combinatorially as well (for large $n$).  Note that (1,1) knots include torus
knots and 2-bridge knots as a proper subset.  In the case of
$(p,pn+1)$ cables, we have the following result:

\begin{theorem}
\label{thm:largen}

Let $K,d$ be as above. Then $\exists \ N>0$ and $c(c',n,p)$ such that $\forall \ n>N$, the following holds:  
$$ {\rm deg}~\HFKa(K_{p,pn+1})=pd+\frac{(p-1)(pn)}{2}. $$
If $i > c(c',n,p)$ we have
$$ \HFKa_*(K_{p,pn+1},i)\cong 
\left\{\begin{array}{ll}

 H_{*+2(k-d)}(\Filt(K,k-d)) & {\text{for
 $i=pd+\frac{(p-1)(pn)}{2}-pk$}} \\

 H_{*+2(k-d)+1}(\Filt(K,k-d)) & {\text{$i=pd+\frac{(p-1)(pn)}{2}-pk-1$}} \\

0 & {\text{otherwise.}}\\
\end{array}
\right. $$
Where $c'$ is a fixed constant coming from the projection of $K$, and
$c(c',n,p)$ is linear in $n$ and quadratic in $p$.
\end{theorem}

In some examples we don't know $H_*(\Filt(K))$. The theorem can still provide useful
information if all that is known is $\HFKa(K)$ in the top filtration dimension.  

\begin{cor}
\label{cor:topfilt}
With $K,d,n>N$ as above 
\begin{align*}
\mathrm{deg}~\HFKa(K_{p,pn+1})&=p\, 
\mathrm{deg}~\HFKa(K) +  \mathrm{deg}~\HFKa(T_{p,pn+1})\\
&=pd+\frac{(p-1)(pn)}{2}.\end{align*}
Furthermore,
\begin{align*}
\HFKa_{*}(K_{p,pn+1},pd+\frac{(p-1)(pn)}{2}) &\cong 
\HFKa_{*-1}(K_{p,pn+1},pd+\frac{(p-1)(pn)}{2}-1)\\&\cong 
\HFKa_{*}(K,d). \end{align*}
\end{cor} 

Of course the corollary is just the restriction of Theorem \ref{thm:largen} to
the top 2 filtration dimensions.   However, it shows that when $\HFKa$
is successful in distinguishing knots by using only the top filtration
dimension (as is the case for the Kinoshita-Terasaka knots and their
Conway mutants, see \cite{GBmutant}, \cite{KinoshitaTerasaka}), it
also distinguishes their $(p,pn+1)$ cables.  In light of
\cite{genusbounds}, the corollary also shows that in many cases the 
Seifert genus of cabled knots is a linear function of the
companion knot's genus, a result proved in \cite{Shibuya} in general.

The proof of Theorem \ref{thm:largen} relies on a special choice of Heegaard diagram for the
cables of a knot which greatly 
simplifies their chain complexes.  This diagram will be introduced in
Section \ref{sec:Heegaards} and will subsequently be used to calculate
$\HFKa(T_{p,q})$ for some of the torus knots.  With the aid of the
diagram and the torus knot calculations, Theorems \ref{thm:2ncable} and
\ref{thm:largen} will be proved in Section \ref{sec:largen}.
Section \ref{sec:examples} will then apply Theorem \ref{thm:2ncable}
to calculate $\HFKa$ for $(2,2n+1)$ cables of the $(2,2m+1)$ torus knots.

\medskip
\noindent{\bf{Remarks}}\qua It is interesting to compare the theorems and corollary above with Equation \eqref{eq:satpoly}. We also remark that all of the above results have corresponding
analogues when $n< 0$, and hence we obtain results for 
$(p,pn\pm1)$ cables.  These are discussed at the end of
Section \ref{sec:largen}. 

\medskip
\noindent{\bf{Acknowledgment}}\qua I cannot thank Peter
Ozsv{\'a}th enough for his willingness and patience to teach me the subject and his
 support and enthusiasm as my advisor.

\section{Preliminaries, Heegaard diagrams, and useful Examples}
\label{sec:Heegaards}

\subsection{Preliminaries on knot Floer homology}
Let $K\subset S^3$ be a knot.  In \cite{Knots}, Ozsv{\'a}th and Sz{\'a}bo introduced the knot
Floer complex $\CFKa(K)=\bigoplus_{i\in\Z}\CFKa(K,i)$ 
associated to a Heegaard diagram for a knot, and whose homology groups
are knot invariants, see
also \cite{RasmussenThesis},\cite{AltKnots}. This complex depends upon a suitable
choice of Heegaard diagram, compatible with the knot in the following
sense:

\begin{defn}
\label{def:hd}
A {\em compatible doubly-pointed Heegaard diagram} for a knot K (or simply a Heegaard diagram for K) is a collection of data 
$$(\Sigma,\{\alpha_1,\ldots,\alpha_g\},\{\beta_1,\ldots,\beta_{g-1},\mu\},w,z),$$
where  
\begin{itemize}
\item $\Sigma$ is an oriented surface of genus g
\item $\{\alpha_1,\ldots,\alpha_g\}$ are pairwise disjoint, linearly
  independent embedded circles which specify a handlebody, $U_\alpha$,
  bounded by $\Sigma$
\item $\{\beta_1,\ldots,\beta_{g-1},\mu\}$ are pairwise disjoint, linearly
  independent embedded circles which specify a handlebody, $U_\beta$,
  bounded by $\Sigma$ such that $U_\alpha\cup_{\Sigma}U_\beta$ is
  diffeomorphic to $S^3$
\item If we do not attach the handle specified by $\mu$, together with
  the final three-ball necessary to make $S^3$, then the resulting
  three-manifold with boundary is the knot complement, $S^3\setminus \nu (K)$
  (i.e.\ $\mu$ is the meridian of the knot)
\item The points $z$ and $w$ can be joined by a small arc $\delta$
 , oriented from $z$ to $w$, which intersects none of
 $\{\alpha_1,\ldots,\alpha_g,\beta_1,\ldots,\beta_{g-1}\}$ and
 algebraically intersects $\mu$ the same number as $lk(K,\mu)$ if we
 arbitrarily orient $\mu$.

\end{itemize}
\end{defn}

We now briefly recall the definitions of $\CFKa(K)$ and its boundary operator in
terms of this diagram, though the reader unfamiliar with the subject
is strongly encouraged to read Section 2 of \cite{AltKnots}.  While
 \cite{Knots} sets up the machinery for knot Floer homology in a more
general context, the level of generality here will be consistent with
that of \cite{AltKnots}.  For this reason all definitions and
notation used here are consistent with those of \cite{AltKnots}
unless otherwise specified.  

Recall that the knot Floer
homology is a doubly graded homology theory.  One grading is a homological grading (coming from the grading on the Floer homology
of $S^3$), also called the Maslov grading.  The other
grading comes from a filtration of $\CFa(S^3)$ induced by the knot.  This latter grading will be referred to
as the filtration or Spin$^c$ grading.  The chain complex $\CFKa$ is
generated by intersection points of the tori
$\T_{\alpha}= \alpha_1 \times \ldots \times \alpha_g$ and
$\T_{\beta}= \beta_1 \times \ldots \times \beta_{g-1} \times \mu$ contained in
Sym$^g(\Sigma_g)$.  Any two such points can be connected by a Whitney
disk $\phi$ whose boundary is contained in the tori. We denote the
intersection number of $\phi$ with the submanifold
$p\times$Sym$^{g-1}(\Sigma_g)$ by $n_p(\phi)$, where $p$ is any point
in $\Sigma-\alpha_1-\ldots-\alpha_g-\beta_1- \ldots-\beta_{g-1}-\mu$.

The relative Maslov and Spin$^c$ gradings (denoted $\gr$
and $\Filt$ respectively) are
determined by the following (found in \cite{AltKnots}):
\begin{gather}
\label{eq:gr}
\gr(\x)-\gr(\y)=\Mas(\phi)-2n_\BasePt(\phi)\\
\label{eq:filt}
\Filt(\x)-\Filt(\y)=n_\FiltPt(\phi)-n_{\BasePt}(\phi),
\end{gather}
\noindent where $\mu$ is the Maslov index of $\phi$.  The
absolute Maslov grading is obtained by the convention that $\gr(\x)=0$
for $\x$ generating $\HFa(S^3)\cong\Z$.  The absolute filtration grading can
be naively obtained by requiring $\HFKa(K,i)$ to be symmetric about
$i=0$, though it has a more invariant description given in
 \cite{Knots}.

 The boundary operator $\partial_z$ is defined as follows: 
$$\partial_z[\x]=\sum_{\y\in\Ta\cap\Tb}
\sum_{\{\phi\in\pi_2(\x,\y)\big|\Mas(\phi)=1,n_w(\phi)=0\}}
\#\left(\UnparModFlow(\phi)\right)[\y],$$
where here $\UnparModFlow(\phi)$ denotes the quotient of the moduli
space of $J$-holomorphic disks representing the homotopy type of
$\phi$, $\ModSp(\phi)$, divided out by the natural action of $\R$
on this moduli space.  This operator can act on various chain
complexes: $\CFa(S^3)$, $\Filt(K,j)$, $\CFKa(K,j)$, for example.  The resulting
homologies will be denoted $\HFa(S^3)$, $H(\Filt(K,j))$, and
$\HFKa(K,j)$.  The operator decomposes as a sum
$\partial_z=\partial_z^0 + \partial_z^1 + \ldots + \partial_z^k$ where
$\partial_z^i$ has the same formula as above except we require
$n_z(\phi)=i$ in addition to $n_w(\phi)=0$. We define $\partial_w$ with the same formula except we
require $n_z(\phi)=0$.  Similarly $\partial_w^i$ requires
$n_z(\phi)=0$ and $n_w(\phi)=i$.

\subsection{A Heegaard diagram for cables}

The purpose of this section is to demonstrate an appropriate Heegaard
diagram for the $(p,pn+1)$ cable of a knot $K$.  The following lemma describes a procedure for finding a Heegaard diagram for cable knots starting from a diagram for the pattern.  Note that, strictly speaking, the diagram we obtain for the cable is not actually compatible in the sense of Definition \ref{def:hd} since it does not contain the meridian of the cable as a $\beta$ attaching curve.  However, a compatible diagram can easily be obtained by stabilizing the diagram described below.  Due to the independence of Heegaard Floer homology under (de)stabilization \cite{HolDisk}, we simply work with the diagram for the cable obtained below.  For more details on the Heegaard diagrams used in this paper see Chapter 2 of \cite{mythesis}.

\begin{figure}[ht!]\small\anchor{fig:unknot}
\begin{center}
\psfrag{a}{\small$\alpha$}
\psfrag{m}{\small$\mu$}
\psfrag{l}{\small$\lambda$}
\psfrag{d}{\small$\delta$}
\psfraga <1.5pt,-1.5pt> {x0}{\small$\x_0$}
\psfraga <0pt,-1pt> {w}{\small$\w$}
\psfraga <-1pt,-1pt> {z}{\small$\z$}
\includegraphics{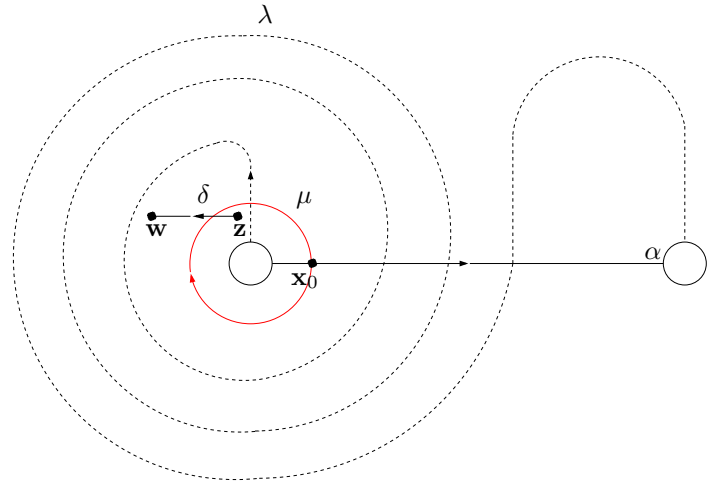}
\end{center}
\caption{\label{fig:unknot} \small Heegaard diagram for the unknot. $\lambda$ represents a
3-framed longitude for the unknot around which we will wind the meridian, $\mu$.}
\end{figure}

\begin{figure}[ht!]\small\anchor{fig:unknotdouble}
\begin{center}
\psfrag{a}{\small$\alpha$}
\psfrag{m}{\small$\tilde{\beta}$}
\psfrag{d}{\small$\delta$}
\psfrag{x0}{\small$\x_0$}
\psfrag{x1}{\small$\x_1$}
\psfrag{x2}{\small$\x_2$}
\psfrag{x3}{\small$\x_3$}
\psfrag{x4}{\small$\x_4$}
\psfrag{x5}{\small$\x_5$}
\psfrag{x6}{\small$\x_6$}
\psfrag{w}{\small$\w$}
\psfrag{z}{\small$\z'$}
\includegraphics[width=0.97\hsize]{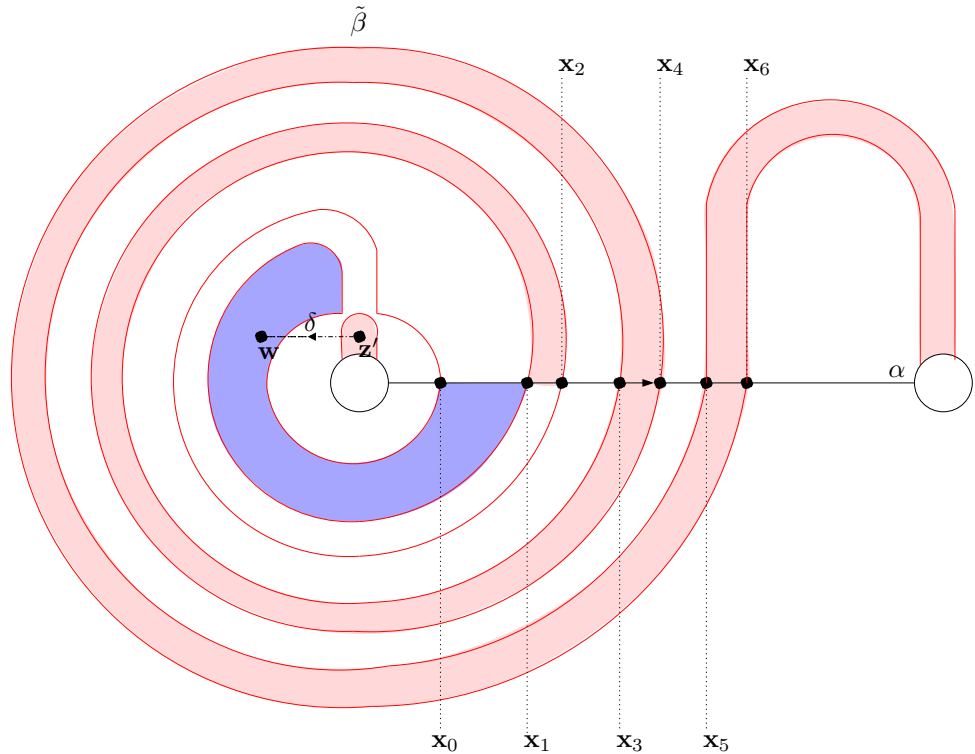}
\end{center}
\caption{\label{fig:unknotdouble}
\small Heegaard diagram for the (2,7) cable of the unknot (i.e.\
  $T_{2,7}$) obtained from the diagram of the unknot
  above via Lemma \ref{lemma:hd}. $\delta$ has intersection 2 with
  $\tilde{\beta}$. The darkly shaded (blue) region indicates the domain
  of the differential $\phi$ connecting $\x_1$ to $\x_0$ discussed in
  the proof of Proposition \ref{prop:T2n} which has
  $n_{w}(\phi)=1$.  The lightly shaded (red) region indicates the domain of
  the differential $\psi$ connecting $\x_1$ to $\x_2$ and which has $n_{z'}(\psi)=1$}
\end{figure}

\begin{lemma}
\label{lemma:hd}
Let
$(\Sigma,\{\alpha_1,\ldots,\alpha_g\},\{\beta_1,\ldots,\beta_{g-1},\mu\},z,w),$
be a Heegaard diagram for a knot K.  Then a Heegaard diagram for
$K_{p,pn+1}$ is obtained from the diagram for K by replacing $\mu$
with a curve $\tilde{\beta}$.  The curve $\tilde{\beta}$ is obtained by winding
$\mu$ along an $n$-framed longitude for the knot $(p-1)$ times. The point $w$
is to remain fixed under this operation. The point $z$ is replaced by a
basepoint $z'$ so that the arc
connecting $z'$ and $w$ has algebraic intersection number $p$
with $\tilde{\beta}$ and is disjoint from all other $\beta$ curves and all $\alpha$ curves.  (See
Figures \fref{fig:unknot} and \fref{fig:unknotdouble}.)
\end{lemma}

\begin{proof}
Let us first understand the Heegaard diagram for $K_{2,1}$ in terms of
the Heegaard diagram for K.  Begin with the unknot.  A Heegaard diagram
for the unknot is simply the standard Heegaard diagram for $S^3$,
together with two points, $z$ and $w$, placed a small distance apart on
either side of the curve $\mu$ representing the meridian, 
\figref{fig:model}A.

\begin{figure}[ht!]\small\anchor{fig:model}
\begin{center}
\psfrag{a}{\small$\alpha$}
\psfrag{m}{\small$\mu$}
\psfrag{d}{\small$\delta$}
\psfrag{z}{\small$\z$}
\psfrag{w}{\small$\w$}
\psfrag{a1}{\small$\alpha_1$}
\psfrag{a2}{\small$\alpha_2$}
\psfrag{z'}{\small$\z'$}
\psfraga <-2pt,0pt> {m'}{\small$\mu'$}
\psfrag{b}{\small$\beta$}
\psfrag{b'}{\small$\tilde{\beta}$}
\psfrag{A}{A}
\psfrag{B}{B}
\psfrag{C}{C}
\psfrag{D}{D}
\psfrag{a'}{\small$\alpha'$}
\psfrag{k}{\small$K$}
\psfrag{k'}{\small$K_{2,1}$}
\includegraphics[width=0.97\hsize]{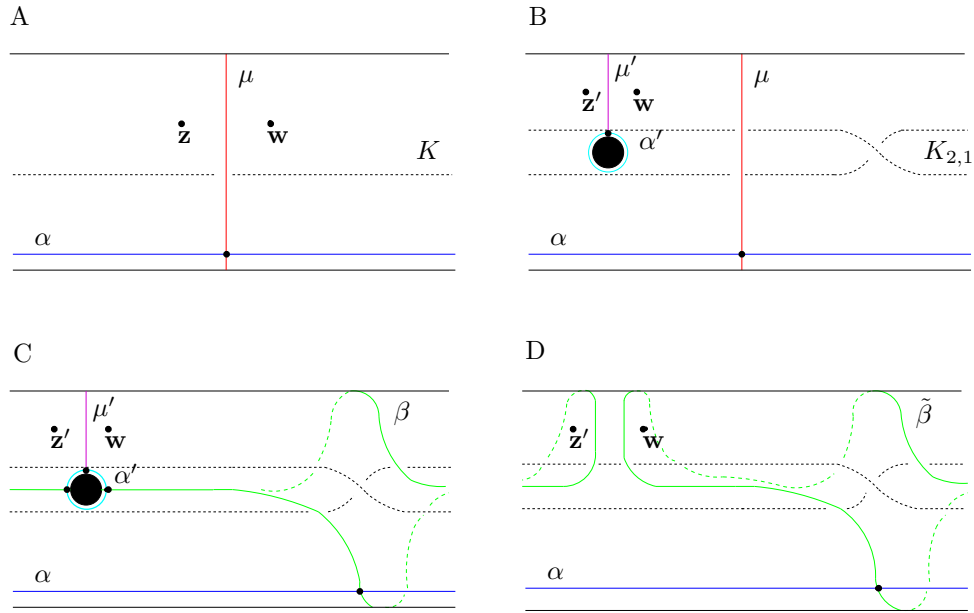}
\caption{\label{fig:model}
Illustration of Lemma \ref{lemma:hd}. The four figures
represent the steps of the lemma.  Each is a top view of
a region of the solid torus. A dashed line represents either the knot
in the core of the solid torus, or an attaching circle for the
Heegaard diagram which is on the underside of the torus.  The dark
circle in figures B and C is a hole drilled through
the torus.  It is destabilized in D after two handleslides.}
\end{center}
\end{figure}

Now draw in place of the unknot its (2,1) cable. Of course this is
still the unknot.  However, $\mu$ is no longer a meridian for the
knot. It has $lk(K,\mu)=2$.  Thus we stabilize the diagram in the sense of
\cite{Gompf} by drilling a hole between the strands of the
knot.  When stabilizing, we add two curves to the diagram,
$\alpha^\prime$ and $\mu^\prime$, which
bound disks when their corresponding handles are attached, and which
satisfy $\alpha^\prime\cap\mu^\prime=1$.  The curve $\mu^\prime$(which does not encircle the added hole) can be chosen so that it is a
meridian for the cabled unknot. On either side of this meridian we add
the points $z^\prime$ and $w$ (we add the prime here and throughout to signify that
this is a diagram for the cabled knot) in such a way to be compatible with the
orientation of the knot.  See \figref{fig:model}B.  Although the diagram with no modification to
the original $\alpha$ and $\mu$ curves represents $S^3$, it is
not a Heegaard diagram for the knot since $\mu$ links the 
knot twice.  Thus we replace $\mu$ by a curve $\beta$ whose attaching
disk does not intersect the knot, and which still results in a
Heegaard diagram for $S^3$.  The attaching disk of $\beta$ stays
between the two strands of the knot, twisting as the knot twists,
while it winds along the longitude of the original unknot. See
\figref{fig:model}C.  Now we perform two handleslides on
$\beta$ to obtain a curve $\tilde{\beta}$ which satisfies the requirements
of the lemma and no longer intersects $\alpha^\prime$. 
We can destabilize the resulting diagram, removing  $\mu^\prime$ and $\alpha^\prime$,to obtain a genus one diagram for the
$(2,1)$ cable of the unknot satisfying the conditions of the
lemma.  See \figref{fig:model}D. 

To find a Heegaard diagram for the $(2,2n+1)$ cable of an arbitrary knot
we note that our choice of the unknot above was not special:
the same sequence of Heegaard moves could have been applied with an arbitrary knot
in a handlebody of genus $g$. The index, $n$, of the
cable corresponds to the framing of the knot used in the Heegaard
diagram.  If we use the zero framing, the procedure yields the $(2,1)$
cable as above.  Picking a framing which adds $n$ meridians to the
$0$-framed longitude for $K$ is equivalent to cabling with the $(2,2n+1)$
torus knot.

Finally, to obtain a diagram for the $(p,pn+1)$ cable of a knot, it is
easy to see (if somewhat harder to draw) that the above procedure can
be extended. The only difference is that instead of winding the
original meridian once along the longitude, we wind it $(p-1)$ times
(while still winding $n$ times in the meridional direction.)
\end{proof}

\subsection{Examples: $\HFKa(T_{2,2n+1}),\HFKa(T_{3,7})$}
\label{subsec:T2n}
As both an illustration of Lemma \ref{lemma:hd} and also as a tool
for the general case, we calculate $\HFKa$ for the $(2,2n+1)$ and $(3,7)$ torus
knots. As it turns out, much of what we see in these simple examples is
reflected in the cables of an arbitrary knot when $n$ is
large.  

\begin{prop}
\label{prop:T2n}
$$\HFKa(T_{2,2n+1},i)\cong 
\left\{\begin{array}{ll}
\Z_{(i-n)} & {\text{for  $|i| \le n$ if $n\ge 0$}} \\
\Z_{(i-n-1)} & {\text{for $|i| \le -n-1$ if $n<0$.}}
\end{array}
\right.$$
The groups in both cases vanish outside of the specified range for i.
\end{prop}

Of course this result is known, \cite{HolDiskTwo}.  (For
generalizations, see
Theorem 1.3 of \cite{AltKnots}, Theorem 1.2 of \cite{Lens}, or \cite{RasmussenThesis}.)  Nonetheless it will be
interesting to see the result arise in this context.

\begin{proof}
To obtain a Heegaard diagram for $T_{2,2n+1}$, we start with the
standard genus one diagram
for the unknot and apply Lemma \ref{lemma:hd} using an $n$-framed
longitude.  This process is pictured in Figures
\fref{fig:unknot} and \fref{fig:unknotdouble}.   From \figref{fig:unknotdouble} we
see there are $2n+1$ intersection points.  To determine the intersection
point which has Maslov grading zero, we disregard the 
basepoint $z'$ in the diagram.  Now we are free to isotope
$\tilde{\beta}$ back around the longitude, removing all
intersection points but $\x_0$.  Thus $\x_0$ generates
$\HFa(S^3)\cong\Z$ and has Maslov grading zero.

We claim that for $n>0$,
$\gr(\x_{i})-\gr(\x_{i+1})=\Filt(\x_{i})-\Filt(\x_{i+1})=1$.  This
follows immediately from the fact that for $i$ odd there is a unique
holomorphic disk $\phi$ connecting $\x_i$ to $\x_{i-1}$ which has
$n_w(\phi)=1$ and $n_{z'}(\phi)=0$, and also a unique holomorphic disk
$\psi$ connecting $\x_i$ to $\x_{i+1}$ which has $n_w(\psi)=0$ and
$n_{z'}(\psi)=1$.  Examples of these disks are shown in \figref{fig:unknotdouble}.  The
statement about the Maslov and filtration gradings follows from Equations
\eqref{eq:gr} and \eqref{eq:filt} above.

For $n\ge 0$, the proposition follows immediately by noting that each intersection point lives in a distinct filtration level.  The case for $n<0$
is completely analogous, except that  
$\gr(\x_{i})-\gr(\x_{i+1})=\Filt(\x_{i})-\Filt(\x_{i+1})=-1$ since the
winding occurs in the opposite direction. Be careful to note that
using a $-1$ framed longitude to obtain a Heegaard diagram for
$T_{2,-1}$ actually produces a Heegaard diagram with 3 intersection
points (when one would expect it to have 1, since it is the
unknot). However, these additional intersection points can be removed
using an allowed isotopy of $\tilde{\beta}$ (i.e.\ one that does not
cross $z'$ or $w$.) Indeed, whenever $n<0$ we can remove two
intersection points in this way. 
\end{proof}

\begin{figure}[ht!]\small\anchor{fig:T37}
\begin{center}
\psfrag{a}{\small$\alpha$}
\psfrag{b}{\small$\tilde{\beta}$}
\psfrag{x0}{\small$\x_0$}
\psfrag{x1}{\small$\x_1$}
\psfrag{x2}{\small$\x_2$}
\psfrag{x3}{\small$\x_3$}
\psfrag{x4}{\small$\x_4$}
\psfrag{x5}{\small$\x_5$}
\psfrag{x6}{\small$\x_6$}
\psfrag{x7}{\small$\x_7$}
\psfrag{x8}{\small$\x_8$}
\psfrag{w}{\small$\w$}
\psfrag{z'}{\small$\z'$}
\includegraphics[width=0.97\hsize]{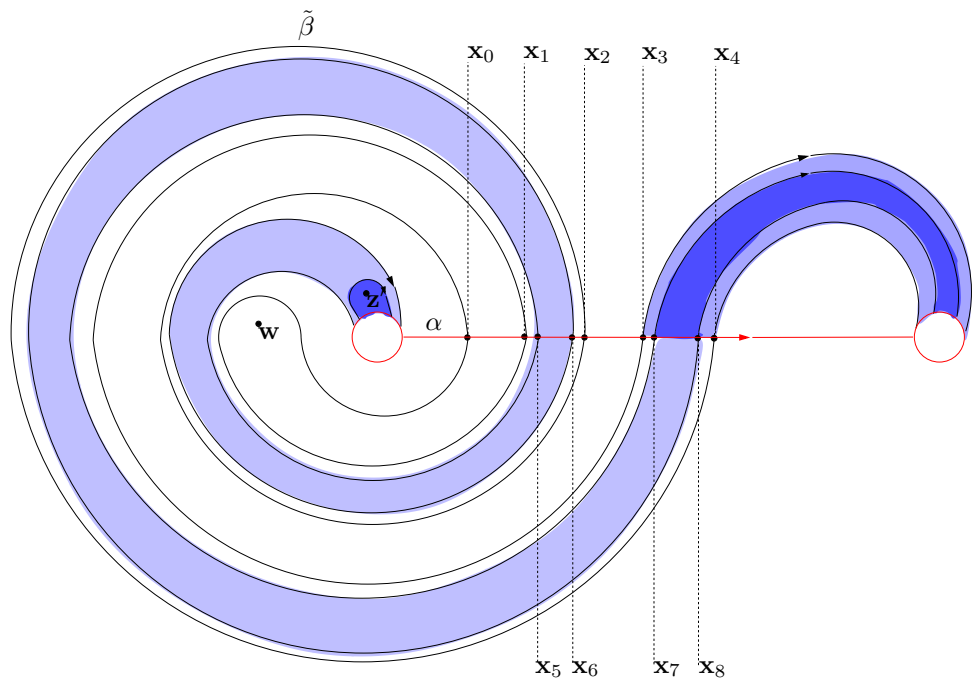}
\end{center}
\caption{\label{fig:T37}\small Heegaard diagram for the (3,7) cable of
  the unknot (i.e.\ $T_{3,7}$) obtained via Lemma \ref{lemma:hd} with $p=3,n=2$ The
  shaded region indicates the domain  $\phi$ connecting $\x_3$ to $\x_4$ discussed in
  the proof of Proposition \ref{prop:T37}.  The light shading indicates multiplicity
  $1$ while the dark shading multiplicity $2$}
\end{figure}

In the interest of being concrete, rather than extend the above example to the case of the torus knots
$T_{p,pn+1}$, we calculate only a specific example, $T_{3,7}$.  The case of
$T_{p,pn+1}$ follows in exactly the same spirit as $T_{3,7}$ and is
only notationally more difficult. In general, $\HFKa(T_{p,q})$
follows from Theorem 1.2 of \cite{Lens}.
\begin{prop}
\label{prop:T37}
$$\HFKa(T_{3,7},i)\cong 
\left\{\begin{array}{ll}

\Z_{(0)} & i=6\\
\Z_{(-1)} & i=5\\
\Z_{(-2)} & i=3\\
\Z_{(-3)} & i=2\\
\Z_{(-4)} & i=0\\
0 & i=1,4, i>6\\
\end{array}
\right.$$
\end{prop}

{\bf Remark}\qua The knot Floer homology groups enjoy a
  symmetry under the natural involution on Spin$^c$ structures given
    by conjugation, Proposition 3.10 of \cite{Knots}. For the
    purposes of this paper, this symmetry can be expressed by the
    following formula:
\begin{equation}
\label{eq:jsym}
\HFKa_{*}(K,i)\cong \HFKa_{*-2i}(K,-i).
\end{equation}
Thus the information in the above proposition completely
specifies $\HFKa(T_{3,7})$.

\begin{proof}
We apply Lemma \ref{lemma:hd} with $p=3$, $n=2$, to obtain the diagram
for $T_{3,7}$ shown in \figref{fig:T37}.  There are nine
intersection points $\x_0 \ldots \x_8$.  Just as in the proof of
Proposition \ref{prop:T2n}, there are unique holomorphic disks with
domains $\phi$ having $n_w(\phi)=1,n_{z'}(\phi)=0$ which connect $\x_1$ to $\x_0$ and
$\x_3$ to $\x_2$.  Both of these domains have $\mu(\phi)=1$.   There are
obvious holomorphic disks connecting $\x_1$ to $\x_2$ and $\x_3$ to $\x_4$ with
$n_{z'}(\phi)=2,n_w(\phi)=0$ and $\mu(\phi)=1$ (see \figref{fig:T37}).  The last two obvious
holomorphic disks connect  $\x_5$ to $\x_6$ and $\x_7$ to $\x_8$ with
$n_{z'}(\phi)=1,n_w(\phi)=0$, and $\mu(\phi)=1$.  

We now calculate the relative filtration difference between $\x_4$ and
$\x_6$.  Let $A$ and $B$ be a symplectic basis for $H_1(\Sigma_1,\Z)$ such
that $[B]\cm[A]=+1$ and so 
that $A=\alpha$ and $[B]=[\tilde{\beta}]$ (assuming
we orient the curves as in \figref{fig:T37}).  Draw an arc from $\x_4$ to $\x_6$ along
$\alpha$ and an arc from $\x_6$ to $\x_4$ following
$\tilde{\beta}$.  The result is a closed curve $\gamma$ which can be chosen 
so that $[\gamma]\cm[A]=1$ and $[\gamma]\cm[B]=1$.   Hence, $$[\gamma-B-A]=[\gamma-\tilde{\beta}-\alpha]=0 \in H_1(\Sigma_1,\Z).$$  Thus there is a
null-homology $\phi$ for the curve $\gamma-\tilde{\beta}-A$ which has
$n_w(\phi)=0$. $n_{z'}(\phi)$ is then the algebraic
intersection number of an arc $\delta$ connecting $z'$ to $w$, with
$\gamma-\tilde{\beta}-A$.  Equivalently, it is the multiplicity of the null-homology constructed at the point $z^{\prime}$. We see from this that $n_{z'}(\phi)=3$
and hence that $\Filt(\x_4)-\Filt(\x_6)=n_{z'}(\phi)-n_w(\phi)=3-0=3$.
A similar analysis shows that $\Filt(\x_5)-\Filt(\x_7)=3$.  Together
with the statements above, we see that each point is in a different filtration dimension and hence generates a $\Z$ summand in the associated graded homology.
Requiring rk$\HFKa(i)$ to be symmetric about $i=0$ yields the statement of the groups in the proposition. 

To find the Maslov grading of the groups we argue similarly to
Proposition \ref{prop:T2n}.  If we forget the reference point $z'$
and slide $\tilde{\beta}$ back, we find that $\x_0$ must have absolute
grading $0$. The theorem then follows from Equations \eqref{eq:gr} and
\eqref{eq:jsym}  together with the Maslov indices of the disks in the first paragraph. 
\end{proof}

\section{Proof of theorems}
\label{sec:largen}

In this section we prove the theorems stated in the introduction.  The
idea behind these theorems is the following: when we perform the operation of
Lemma \ref{lemma:hd}, the resulting Heegaard diagram can be
simultaneously viewed as a diagram for both the original and
cabled knot by appropriately placing three basepoints.  When we
increase $n$  we add a large number of spirals to the Heegaard
diagram. If we make this spiraling region large enough, we can
ensure that all generators in the complex for both the original and
cabled knots having high filtration gradings live in
the spiraling region.  Since the domains of the differentials in the Heegaard
diagram are the same regardless of whether we are viewing it as a
diagram for the original or cabled knot, we can use our assumed
knowledge of the original knot's differential to calculate the
differential for the cabled knot. When we specialize to the case of
$(2,2n+1)$ cables, the symmetry of $\HFKa$ under the
conjugation action on Spin$^c$ structures (Equation \eqref{eq:jsym}) allows us to completely determine $\HFKa(K_{2,2n+1})$. With the idea in place, we begin.

\subsection{Proof of Theorems \ref{thm:largen} and \ref{thm:2ncable}}

\begin{figure}[ht!]\small\anchor{fig:trefoil}
\begin{center}
\psfrag{m}{\small$\mu$}
\psfrag{d}{\small$\delta$}
\psfraga <-2pt,0pt> {x0}{\small$x_0$}
\psfrag{z}{\small$\z$}
\psfraga <-1pt,-1pt> {w}{\small$\w$}
\psfrag{a1}{\small$\alpha_1$}
\psfrag{a2}{\small$\alpha_2$}
\psfrag{b1}{\small$\beta_1$}
\psfraga <-2pt,0pt> {w1}{\small$w_1$}
\psfraga <-2pt,0pt> {w2}{\small$w_2$}
\psfraga <-2pt,0pt> {w3}{\small$w_3$}
\psfrag{y1}{\small$y_1$}
\psfrag{y2}{\small$y_2$}
\psfraga <-2pt,0pt> {y3}{\small$y_3$}
\includegraphics{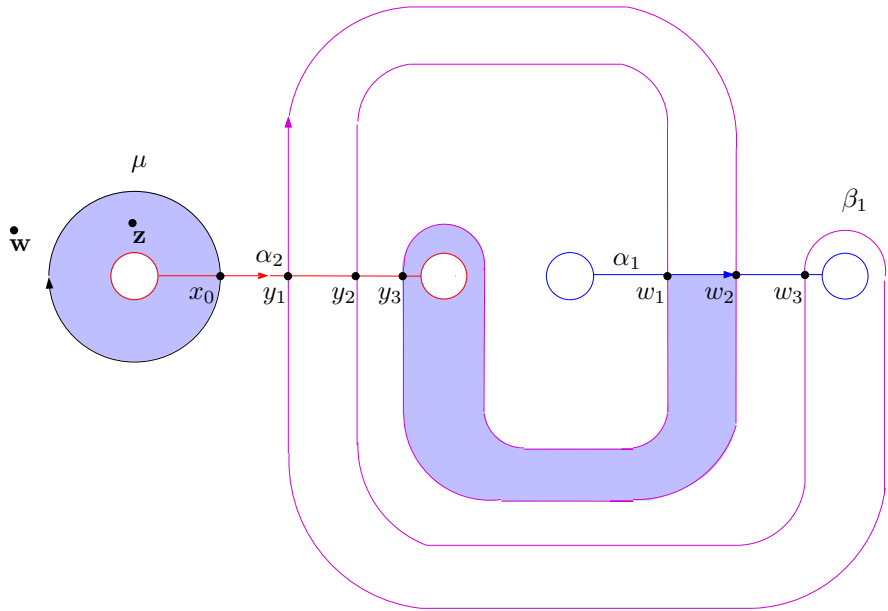}
\caption{\label{fig:trefoil}
Heegaard diagram for the right-handed trefoil coming from its 2-bridge
presentation. Shaded is a domain $\phi$ connecting $(x_0,w_2)$ to
$(x_0,w_1)$ having $n_z(\phi)=1$, $n_w(\phi)=0$.  When we apply Lemma
 \ref{lemma:hd} this domain winds along the longitude and the
resulting domain $\phi'$ has $n_{z'}(\phi')=p$, $n_z(\phi')=1$, and
$n_w(\phi')=0$. See Figure \ref{fig:trefoildouble} for an
illustration when $p=2$.}
\end{center}
\end{figure}

 Given a $g$-bridge presentation
for $K$, we obtain a genus $g$ Heegaard diagram for $S^3$ compatible
with $K$ (see \cite{Gompf}).  See \figref{fig:trefoil} for
the Heegaard diagram of the right-handed trefoil.  Note that the
meridian for the knot will intersect $\alpha_g$ exactly once in
the point $x_0$, and
will intersect none of the other $\alpha$ curves.  This ensures that
all intersection points of the knot's chain complex will be 
of the form $(x_0,\bf\y)$ for some $(g-1)$ tuple, $\bf\y$, of
intersection points.  Following \cite{Knots} one can determine the
filtration grading of these generators.  For each filtration summand,
$\CFKa(K,i)$, there is a unique set of intersection points generating
the summand.  Let us denote by $C(i)$ all $g-1$ tuples of
intersection points between $\T_\beta\setminus\mu$ and $\T_\alpha\setminus\alpha_g$
which, together with $x_0 \in \alpha_g\cap\mu$, generate $\CFKa(K,i)$.
Thus we write  $x_0 \times C(i)$ to mean $\CFKa(K,i)$.

\begin{figure}[ht!]\small\anchor{fig:trefoildouble}
\begin{center}
\psfrag{a1}{\small$\alpha_1$}
\psfrag{m}{\small$\tilde{\beta}$}
\psfrag{d}{\small$\delta$}
\psfrag{x0}{\small$x_0$}
\psfrag{x1}{\small$x_1$}
\psfrag{x2}{\small$x_2$}
\psfrag{x3}{\small$x_3$}
\psfrag{x4}{\small$x_4$}
\psfrag{y}{\small$\z$}
\psfrag{y1}{\small$y_1$}
\psfrag{q}{\small$z_1$}
\psfrag{w1}{\small$w_1$}
\psfrag{w3}{\small$w_3$}
\psfrag{w2}{\small$w_2$}
\psfrag{b1}{\small$\beta_1$}
\psfrag{z'}{\small$\z'$}
\psfrag{w}{\small$\w$}
\psfrag{z}{\small$\z$}
\includegraphics[width=0.97\hsize]{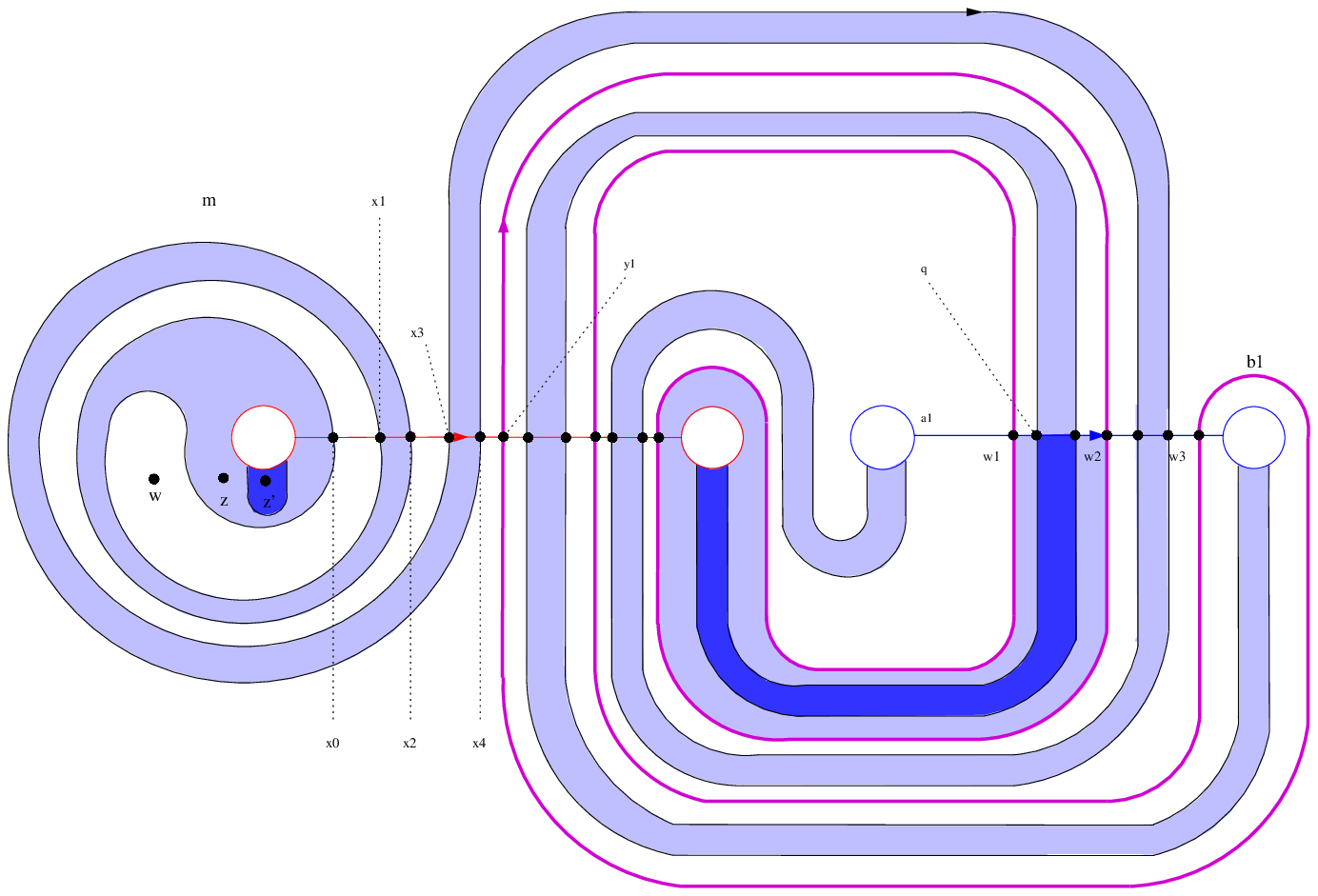}
\caption{\label{fig:trefoildouble} Heegaard diagram for the $(2,13)$ cable of
  the right-handed trefoil. The exterior intersection points are those
  of the form $(x_i,w_j)$ for $i=0,\ldots,4, j=1,2,3$ (there are only
  $5$ points $x_i \in \alpha_2\cap\tilde{\beta}$ generating exterior
  points because the others
  were removed by an isotopy of $\tilde{\beta}$ not crossing $\z'$.)
  Note the domain $\phi'$ from $(x_i,w_2)$ to $(x_i,w_1)$ for
  $i=0,\ldots,4$. This is the domain from Figure \ref{fig:trefoil}
  which was wound along the longitude to have $n_{z'}(\phi')=2$.
  The point $(y_1,z_1)$ is an example of an interior point.}
\end{center}
\end{figure}

Apply Lemma
 \ref{lemma:hd}.  
The first thing to note is that by
appropriately placing a third basepoint, $z$, the
diagram can be viewed as compatible either with $K$ or with
$K_{p,pn+1}$:

\begin{lemma}
Let $z$ be a point on the arc $\delta$ connecting $z^\prime$ to $w$
for which the segment of $\delta$ connecting $z$ to $w$ has
intersection number 1 with $\tilde{\beta}$. See 
 \figref{fig:trefoildouble}.  Then the Heegaard diagram
with the pair $(w,z^\prime)$ is compatible with $K_{p,pn+1}$, while
with the pair $(w,z)$ it is compatible with $K$. 
\end{lemma}
\begin{proof}
That the diagram with $(w,z^\prime)$ is compatible with $K_{p,pn+1}$
is just Lemma \ref{lemma:hd}.  To see that the diagram with $(w,z)$ is
compatible with $K$, simply isotope $\tilde{\beta}$
in the reverse direction to that of Lemma \ref{lemma:hd} in order to arrive
at the original diagram for $K$.  For the diagram with $(w,z)$ this is an
allowed isotopy in the sense of \cite{Knots} since it does not cross
either basepoint.
\end{proof}

Let $n\gg 0$. This creates a large
spiraling region in the Heegaard diagram, similar to the diagrams for
the torus knots in Section \ref{sec:Heegaards}.  This region contains an
odd number, $2(p-1)n+1$, of intersection points of $\tilde{\beta}$ with
$\alpha_g$.  Denote these intersection points $x_0 \ldots
x_{2(p-1)n}$.  See \figref{fig:T37} or \fref{fig:trefoildouble}. (One should be careful here.  If we start with a
  $0$-framed longitude and apply Lemma \ref{lemma:hd}  there will be
  at least $2(p-1)n+1$ intersection points, but possibly more. The
  $0$-framed longitude for $K$ may intersect $\alpha_g$) When we look at $x_i \times C(j)$ we find that
  there is a ``copy'' of each filtration summand of the original knot's chain
  complex carried by the point $x_i$. In light of this we define:

\begin{defn}
In the Heegaard diagram for the $(p,pn+1)$ cable of $K$ we call an
intersection point an {\em exterior intersection point} if it is of the form
$(x_i,\y)$ where $x_i \in \tilde{\beta}\cap \alpha_g$ and $x_i$ can be
joined to $x_0$ by an arc which intersects $\tilde{\beta}$
geometrically at most $2(p-1)n-1$ times and which intersects none of
the other attaching curves.  All other intersection points will be
called {\em interior}.  
\end{defn}

The exterior points are those that arise from the spiraling region.
Before continuing further we establish a
labeling convention for the $\alpha_g \cap\tilde{\beta}$ component of
the exterior points.  We label these points as follows:
each time we wind $\tilde{\beta}$ around the longitude (i.e.\ increase the parameter $p$) we add $2n$ intersection points $x_i \in
\alpha_g\cap \tilde{\beta}$ generating exterior points.  First we label the points that arise the
first time we wind around the longitude. From left to right we label
them $x_0\ldots x_{2n}$. If $p=2$, we are done.  If $p>2$ we next label the points arising the
second time we wind around the longitude $x_{2n+1},\ldots,x_{4n}$, again from left to right,
and so on. See
\figref{fig:T37} for a picture of the spiraling region with this
labeling convention.

The exterior intersection points play a primary role in the
proof of the theorem. If we make $n$ large enough, the chain complexes
in highest
filtration dimensions (with respect to the filtrations induced by both
the uncabled and cabled knot) will be generated by a subset of the
exterior points.   To see this, we must understand the relative
filtrations induced by $K$ and $K_{p,pn+1}$.  We develop this
knowledge through a sequence of lemmas. Let us denote the filtration
with respect to $z$ (i.e.\ induced by $K$) by $\Filt$ and the filtration with respect
to $z'$ (i.e.\ induced by $K_{p,pn+1}$) by $\Filt'$.

We begin with the relative filtration between exterior points sharing the same
$(g-1)$-tuple, $\y$.

\begin{lemma}
For $i < 2n$ odd, we have 
\begin{gather*}
\Filt(x_{i-1} ,\y)-\Filt(x_{i},\y)=\Filt'(x_{i-1},\y)-\Filt'(x_{i},\y)=1\\
\Filt(x_i,\y)-\Filt(x_{i+1},\y)=0\\
\Filt'(x_i,\y)-\Filt'(x_{i+1},\y)=p-1.
\end{gather*}
\end{lemma}

\begin{proof}
For $i$ odd, there is a holomorphic disk with domain $\phi$ from $(x_i,\y)$ to $(x_{i-1},\y)$ having
$n_z(\phi)=n_{z'}(\phi)=0$, $n_w(\phi)=1$. It is the product of the
disk from Propositions \ref{prop:T2n} and \ref{prop:T37} (connecting $\x_i$ to
$\x_{i-1}$) with the constant map in Sym$^{g-1}(\Sigma_g)$.   There is
a disk with domain $\psi$ from
$(x_i,\y)$ to $(x_{i+1},\y)$ with $n_z(\psi)=n_w(\psi)=0$,
$n_{z'}(\psi)=p-1$. This is the the analogue of the disk in
Propositions \ref{prop:T2n} and \ref{prop:T37} connecting $\x_i$ to
$\x_{i+1}$. Topologically it is still a disk, but now it wraps around
the longitude $p-1$ times. The lemma follows from Equation \eqref{eq:filt}.
\end{proof}

Next we fix $x_i$ and vary the $(g-1)$-tuple.
\begin{lemma}
\label{lemma:varytuple}
Suppose $\y \in C(j), \z \in C(k)$.  Then,
\begin{gather*}\Filt(x_i,\y)-\Filt(x_i,\z)=j-k\\
\Filt'(x_i,\y)-\Filt'(x_i,\z)=p(j-k).\end{gather*}
\end{lemma}

\begin{proof}
The two $(g-1)$-tuples had
filtration difference $j-k$ in the original Heegaard diagram
 by assumption (i.e.\ before we applied Lemma \ref{lemma:hd}).  Thus the boundary of the domain connecting
these points in the original diagram had intersection number $j-k$
with an arc connecting $z$ and $w$.  In the cabled diagram, the boundary
of the new domain will still have intersection number $j-k$ with the
arc connecting $z$ to $w$ while it will have intersection number
$p(j-k)$ with the arc connecting $z'$ and $w$. See Figures
 \fref{fig:trefoil} and \fref{fig:trefoildouble} for an illustration of
this lemma.
\end{proof}

\begin{figure}[ht!]\small\anchor{fig:chaincomplex}
\psfraga <0pt,-5pt> {1}{\small $ C(d)$}
\psfraga <0pt,-5pt> {a}{\small $x_0$ }
\psfraga <0pt,-5pt> {2}{\small $ C(d-1)$}
\psfraga <0pt,-5pt> {b}{\small $x_1$}
\psfraga <0pt,-5pt> {3}{\small $C(d-2)$}
\psfraga <0pt,-5pt> {c}{\small $x_2$}
\psfraga <0pt,-5pt> {4}{$\ldots$}
\psfraga <0pt,-5pt> {5}{\small $C(-d)$}
\psfraga <0pt,-5pt> {d}{\small $x_3$}
\psfraga <0pt,-5pt> {e}{\small $x_4$}
\psfraga <0pt,-5pt> {f}{\small $\vdots$}
\psfraga <0pt,-5pt> {g}{\small $x_{2n}$}
\psfraga <0pt,-5pt> {a1}{\small $(0,0)$}
\psfraga <0pt,-5pt> {a2}{\small $(-1,-p)$}
\psfraga <0pt,-5pt> {a3}{\small $(-2,-2p)$}
\psfraga <0pt,-5pt> {a5}{\small $(-2d,-2pd)$}
\psfraga <0pt,-5pt> {b1}{\small $(-1,-1)$}
\psfraga <0pt,-5pt> {b2}{\small $(-2,-p-1)$}
\psfraga <0pt,-5pt> {b3}{\small $(-3,-2p-1)$}
\psfraga <0pt,-5pt> {b5}{\small $(-2d-1,-2pd-1)$}
\psfraga <0pt,-5pt> {c1}{\small $(-1,-p)$}
\psfraga <0pt,-5pt> {c2}{\small $(-2,-2p)$}
\psfraga <0pt,-5pt> {c3}{\small $(-3,-3p)$}
\psfraga <0pt,-5pt> {c5}{\small $(-2d-1,-2pd-p)$}
\psfraga <0pt,-5pt> {d1}{\small $(-2,-p-1)$}
\psfraga <0pt,-5pt> {d2}{\small $(-3,-2p-1)$}
\psfraga <0pt,-5pt> {d3}{\small $(-4,-3p-1)$}
\psfraga <0pt,-5pt> {d5}{\small$(-2d-2,-2pd-p-1)$}
\psfraga <0pt,-5pt> {e1}{\small $(-2,-2p)$}
\psfraga <0pt,-5pt> {e2}{\small $(-3,-3p)$}
\psfraga <0pt,-5pt> {e3}{\small $(-4,-4p)$}
\psfraga <0pt,-5pt> {e5}{\small $(-2d-2,-2pd-2p)$}
\psfraga <0pt,-5pt> {g1}{\small $(-n,-np)$}
\psfraga <0pt,-5pt> {g2}{\small $(-1-n,-np-p)$}
\psfraga <0pt,-5pt> {g3}{\small $(-2-n,-np-2p)$}
\psfraga <0pt,-5pt> {g5}{\small $(-2d-n,-2pd-np)$}
\includegraphics[width=0.97\hsize, height=2.9in]{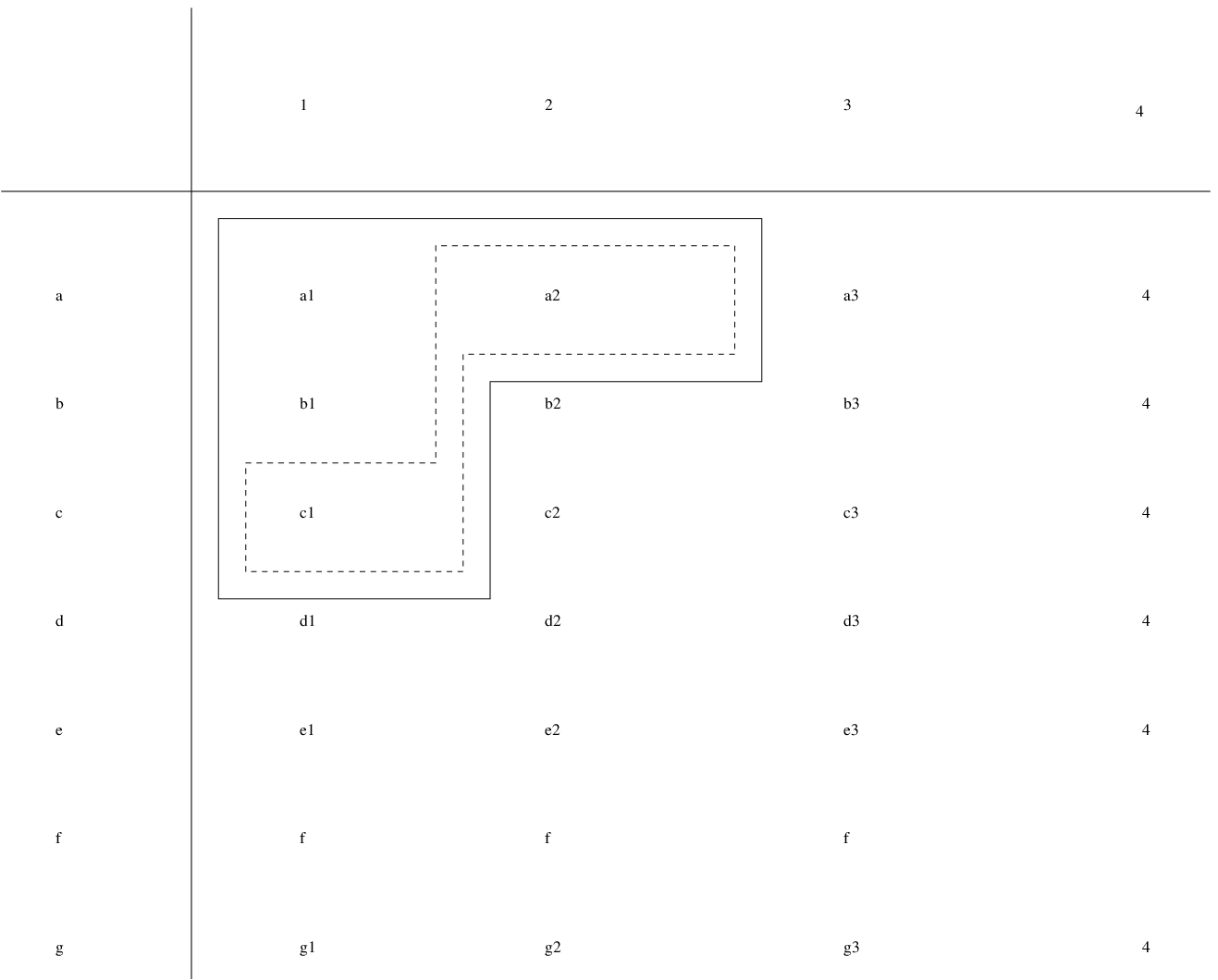}
\caption{\label{fig:chaincomplex}
Table of relative filtrations of exterior points for $K$ and
  $K_{p,pn+1}$ when $n>N$.  The columns and rows are arranged to distinguish the
  elements $(x_i,\y)\in x_i\times C(j)$. The number on the left is the relative filtration of the
  points in each summand taken with respect to $z$ (the original knot's  basepoint).  The number on the right is the relative
  filtration of the points taken with respect to $z'$ (the cabled
  knot's basepoint).  The dashed box is the chain complex for
  the cabled knot in relative filtration $-p$.  It is contained in the
  solid box of generators with relative $z$ filtration greater than
  $-2$.  Proposition \ref{prop:largen2} shows that the complex in the
  dashed box is actually a subcomplex of the complex in the solid box. Note the filtration
  gradings here are relative.
}

\end{figure}

We now know the relative filtration grading for both the cabled and
uncabled knot of all intersection points of the form $(x_i,\y)$, with
$i\le 2n$. See \figref{fig:chaincomplex} for a
table depicting the two chain complexes.

There are many more intersection points in general -- $\tilde{\beta}$ intersects the other $\alpha$ curves as well
as $\alpha_g$. See \figref{fig:trefoildouble} for an example of these intersection
points. We must also understand the filtration grading of the points $(x_i,\y)$
when $i>2n$ (the rest of the exterior points).  These points occur
when $p>2$ because $\tilde{\beta}$ winds more than once around the longitude for $K$. To this end we have the following:

\begin{lemma}
\label{lemma:isolate}
There exist constants  $N,c^\prime > 0$ such that $\forall \ n$ with $n> N$,
the relative filtration gradings of the exterior intersection points
containing $x_0 \ldots x_{2(n-c')}$  are higher than the filtration grading of all other intersection points
except possibly those containing $x_{2(n-c')+1} \ldots x_{2n}$. The $n$, as
always, refers to the parameter specifying the cabled knot
$K_{p,pn+1}$. The constant $c^\prime$ is independent of $n$ and depends only on the projection of $K$.
\end{lemma}

\begin{proof}
First note that the relative filtration difference between any interior
points, $\p$ and $\q$, does not change as we vary $n$.  This is because domains connecting interior points
either remain fixed (take place entirely in the interior of the
diagram) or simply add area to the part of the domain with
multiplicity in the spiral.  In either case $n_{z'}(\phi),\ n_z(\phi), \ n_w(\phi)$
all remain fixed. 

Next observe that the relative filtration difference between points of the form $(x_{2n},\y)$ and
all the interior points is fixed as we vary $n$.  The reason is the same
as above: Domains connecting $(x_{2n},\y)$ to the interior at worst change
by adding area to the part of the domain in the spiral and hence
$n_{z'}(\phi),n_z(\phi),n_w(\phi)$ are all constant.

Finally we must account for the other exterior points. Let $i>0$.
 We need to calculate the filtration difference between $x_{2n}$ and $x_{2n+i}$. 
We claim that both $\Filt(x_{2n})-\Filt(x_{2n+i})$ and
$\Filt'(x_{2n})-\Filt'(x_{2n+i})$ are bounded below by some constant
K, independent of $n$. 

To prove the claim observe that  $\Filt(x_{2n})-\Filt(x_{2n+1})$ and
$\Filt'(x_{2n})-\Filt'(x_{2n+1})$ are independent of $n$ for the same
reasons as above.  The homological method for calculating filtration
differences of Proposition \ref{prop:T37} shows that when $n$ is large enough,
$\Filt(x_{2n+1})-\Filt(x_{2n+j})$ and
  $\Filt'(x_{2n+1})-\Filt'(x_{2n+j})$ are both always greater than or
    equal to zero for any $j>1$.  These observations, the
    previous two lemmas, and the fact that
    there are only a finite number of interior intersection points
    proves the lemma.
\end{proof}

We now wish to prove:

\begin{lemma}
\label{lemma:partial1} For $n>N$ as above,
\begin{gather*}
\HFKa_*(K_{p,pn+1},{\rm topmost}) \cong\HFKa_{*-1}(K_{p,pn+1},
{\rm topmost}-1) \cong \HFKa_*(K,d)\\
\HFKa(K_{p,pn+1},{\rm topmost}-i)\cong0 \quad\text{for}\quad i=2,\ldots,p-1,
\end{gather*}
\noindent Where {\rm topmost} refers to the highest filtration dimension for which
\newline $\HFKa(K_{p,pn+1})\ne0$ and  $d={\rm deg}~\HFKa(K)$.
\end{lemma}

\begin{proof}
When $n>N$ the above Lemmas show that the
exterior points $x_0\times C(d)$ are
higher in relative filtration than all other
intersection points.  This holds whether we take our filtration with
respect to $z$ (the original knot) or $z'$ (its cable). It follows from this and Equation \eqref{eq:filt} that the domain of any disk connecting points of the form $x_0\times C(d)$ with $n_w=0$ must have $n_z=n_{z'}=0$ as well.  This
implies that 
the differential restricted to these points is independent of the basepoint used and we immediately have
$\HFKa(K_{p,pn+1},$topmost$) \cong \HFKa(K,d)$.  We do not yet know the absolute filtration grading
for the cable and hence we simply refer to the dimension as ``topmost''
for now.  The gradings are the same for either knot because the chain
complex calculates $\HFKa(S^3)$ regardless of which way it is
filtered. Thus a generating point for this homology is independent of the filtration.  Furthermore, the relative Maslov grading is calculated using Equation
\eqref{eq:gr}, which is also independent
of the filtration used. This proves that the first and last groups
stated in the lemma are isomorphic.

We show the first and second groups are isomorphic.
Recall that $\partial_w$ decomposes as a sum $\partial_w=\partial_w^0 +
\partial_w^1 + \ldots + \partial_w^k$, where $\partial_w^i$ is the
boundary operator counting holomorphic disks whose domains have
$n_z(\phi)=0$ and $n_w(\phi)=i$.  $(\partial_w)^2=0$ implies
$\partial_w^1$ is a chain map from $\CFKa(K_{p,pn+1},$topmost-1$)$ to
$\CFKa(K_{p,pn+1},$topmost$)$. For each (g-1) tuple $\y$, there is an obvious holomorphic disk
with domain $\phi$
  connecting $(x_1,\y)$ to $(x_0,\y)$  with $n_z(\phi)=n_{z'}(\phi)=0,\
  n_w(\phi)=1$.  It
  is the product, $u\times\y$, of the disk $u$ in the torus from Propositions
  \ref{prop:T2n} and \ref{prop:T37},  with the constant map in
  Sym$^{g-1}(\Sigma_g)$.  We denote this summand in $\partial_w^1$ by
  $l_0$.  In the standard way (see, for instance, Theorem 4.1 of
  \cite{Knots}) we filter the chain map $\partial_w^1$ with respect to
  negative area of the domains of disks.  With respect to this
  filtration, 
$$\partial_w^1=l_0 + \text{lower order terms}.$$
$l_0$ is clearly an isomorphism of chain complexes, and hence
$\partial_w^1$  induces an isomorphism of groups.  The grading shift
is a consequence Equation \eqref{eq:gr}.

From Lemma \ref{lemma:isolate} and \figref{fig:chaincomplex} it is clear that
$\HFKa(K_{p,pn+1},{\rm topmost}-i)\cong0$ for $i=2,\ldots,p-1$: there are simply
  no intersection points in these filtration dimensions.
\end{proof}

\begin{prop}
\label{prop:largen2}

Let K be a knot in $S^3$, and suppose that $\mathrm{deg}~\HFKa(K)=d$. Then
$\exists \ N>0$ such that for all 
$n>N$, the following holds:  
$$\mathrm{deg}~\HFKa(K_{p,pn+1})=pd+\frac{(p-1)(pn)}{2} $$
Furthermore, $\exists \ c(c',n,p)$ such that if $i > c(c',n,p)$ we have
$$ \HFKa_*(K_{p,pn+1},i)\cong 
\left\{\begin{array}{ll}

 H_*(C(K,i'\ge d-k),\partial_w) & {\text{for
 $i=pd+\frac{(p-1)(pn)}{2}-pk$}} \\

\HFKa_{*+1}(K_{p,pn+1},i+1) & {\text{$i=pd+\frac{(p-1)(pn)}{2}-pk-1$}} \\

0 & {\text{otherwise}}\\

\end{array}
\right. $$ 
\noindent Where $C(K,i'\ge d-k)$ is the chain complex generated by
points with filtration dimension $\ge d-k$ with respect to $K$.
$c(c',n,p)$ is linear in $n$ and quadratic in $p$.

\end{prop}

\begin{proof}
The base case being established by Lemma \ref{lemma:partial1} we
assume the theorem holds for $k-1$, with $k<n-c'$ (here $c'$ is the constant from Lemma \ref{lemma:isolate}). It is immediate that for $k<n-c'$
$$ \HFKa_*(K_{p,pn+1},\mathrm{topmost}-pk)) \cong \HFKa_{*-1}(K_{p,pn+1},\mathrm{topmost}-pk-1).$$
$\partial_w^1$ induces an isomorphism of these groups just as in
Lemma \ref{lemma:partial1}.  We will establish the
absolute filtration dimension at the end.  We first note that 
$(C(K,i'\ge d-k),\partial_w)$ is a subcomplex of
$(\CFKa(K),\partial_w)$.  This follows from the fact that $\partial_w$
counts only those disks whose domains have $n_z(\phi)=0$.  Thus the domain of any
differential with range outside of $(C(K,i'\ge d-k),\partial_w)$ must have
negative multiplicity by Equation \eqref{eq:filt}. This is impossible
by Lemma 3.2 of \cite{HolDisk}.  

The generators of $\CFKa(K_{p,pn+1},\mathrm{topmost}-pk)$ are
clearly contained in  \linebreak $(C(K,i'\ge d-k),\partial_w)$, see 
 \figref{fig:chaincomplex}.  We  will show
that the former is actually a subcomplex of the latter.  The differential on
$\CFKa(K_{p,pn+1},\mathrm{topmost}-pk)$ counts disks whose domains have
$n_{z'}(\phi)=n_w(\phi)=0$, while \newline $(C(K,i'\ge d-k),\partial_w)$
requires
$n_z(\phi)=0$.  When restricted to points in
$\CFKa(K_{p,pn+1},\mathrm{topmost}-pk)$, however, both differentials
are identical -- the relative filtration gradings of generators in this set are
the same with respect to either basepoint. Thus Equation \eqref{eq:filt}
implies that both differentials compute the same homology,
$\HFKa(K_{p,pn+1},\mathrm{topmost}-pk)$.  

 To show that  $\CFKa(K_{p,pn+1},\mathrm{topmost}-pk)$ is a
subcomplex is to show that there are no differentials from
$$\CFKa(K_{p,pn+1},\mathrm{topmost}-pk)\quad
\text{to}\quad \frac{(C(K,i'\ge d-k),\partial_w)}{
\CFKa(K_{p,pn+1},\mathrm{topmost}-pk)}.$$ This is an application of Equation
\eqref{eq:filt}.  The numbers on the left in \figref{fig:chaincomplex}
depend on Equation \eqref{eq:filt} using $z$, while the numbers on the
right use $z'$.  If
$$\x\in \CFKa(K_{p,pn+1},\mathrm{topmost}-pk)\ \ \mathrm{and}\ \ \y \in \frac{(C(K,i'\ge
  d-k),\partial_w)}{\CFKa(K_{p,pn+1},\mathrm{topmost}-pk)},$$
 $$\Filt(\x)-\Filt(\y) >
\Filt'(\x)-\Filt'(\y).\leqno{\text{then}} $$ This immediately implies the domain connecting
$\x$ to $\y$ has negative multiplicity (since $n_z(\phi)=0$).  Again
Lemma 3.2 of \cite{HolDisk} implies there are no holomorphic
representatives of these domains and hence $\CFKa(K_{p,pn+1},\mathrm{topmost}-pk)$
is a subcomplex.

The proposition, filtration dimensions aside, follows immediately -- the
quotient complex has trivial homology since
$\partial_w^1$ induces an isomorphism between
\linebreak $\HFKa_*(K_{p,pn+1},\mathrm{topmost}-p(k-1))$ and
$\HFKa_*(K_{p,pn+1},\mathrm{topmost}-p(k-1)-1)$.  

Calculating the absolute filtration dimension follows from
what we have proved and Equations \eqref{eq:satpoly} and
 \eqref{eq:Euler}. The knots $K$ and $K_{p,pn+1}$ have isomorphic
groups in top filtration dimension.  If the Euler characteristic of this
group is non-zero, then so is the coefficient of $T^{\Filt(\mathrm{topmost})}$
  in the Alexander polynomial for both knots, and the result
follows from Equation \eqref{eq:satpoly}.  If the Euler characteristic
of the group is zero, proceed to the first group with non-zero Euler
characteristic. One must exist since the Alexander polynomial of a
knot cannot be zero.  Inspection of the relative filtration levels
for both knots shows that if deg$\HFKa(K)-\mathrm{deg}~\Delta_K(t)=l$, then
deg$\HFKa(K_{p,pn+1})-\mathrm{deg}~\Delta_{K_{p,pn+1}}=p\cm l$. This completes
the proposition.
\end{proof}

The constant $c(c',n,p)$ in the proposition is explained as follows: Lemma \ref{lemma:isolate} shows that for $n{>}N$, the exterior
points carried by $x_0, \ldots, x_{2(n-c')}$ are higher in filtration dimension than
all points except those carried by $x_{2(n-c')+1},\ldots,x_{2n}$, with
$c'>0$ independent of $n$.   It follows that 
the exterior points $(x_i,\y)$ with $i<2n$ generate the top $p(n-c^\prime)+1$
filtration dimensions.  Since
deg$\HFKa(K_{p,pn+1})=pd+\frac{(p-1)(pn)}{2}$ we see that the constant
$c(c',n,p)=pd+\frac{(p-1)(pn)}{2}-p(n-c')-1$ takes the appropriate form. 

\begin{proof}[Proof of Theorems \ref{thm:2ncable} and
  \ref{thm:largen}] 
Theorem \ref{thm:largen} is a restatement of Proposition
 \ref{prop:largen2}.  The chain complex $(C(K,i'\ge
d-k),\partial_w)$ is naturally identified with $\Filt(-K,k-d)$, where
$-K$ denotes $K$ with the reverse orientation -- the differential on
$\Filt(-K)$ is $\partial_w$ by definition, while the relative
filtration equation (Equation \eqref{eq:filt}) for $-K$ permutes $z$
and $w$, see \cite{Knots}.  The homology of $\Filt(-K,k-d)$ is identified with
$\Filt(K,k-d)$ in Section 3.5 of \cite{Knots}.  

The grading shift occurs for the following reason: the
relative gradings for both the original and cabled knots are defined
using 
$$\gr(\x)-\gr(\y)=\mu(\phi)-2n_w(\phi).$$
The proof of Proposition \ref{prop:largen2}, however, identifies
the knot Floer homology groups of the cable with the  homology of the complexes, $(C(K,i'\ge
d-k),\partial_w)$.  Using the differential, $\partial_w$,
the relative grading equation is
$$\gr(\x)-\gr(\y)=\mu(\phi)-2n_z(\phi).$$
The grading shift is a consequence of this and the
relative filtrations.

In the special case
$p=2$, we see that the constant $c(c',n,p)$ is negative for
sufficiently large $n$, thus proving Theorem \ref{thm:2ncable}.
\end{proof}

{\bf Remark}\qua  In the case where $n<0$, the above discussion
carries through almost verbatim.  The exterior points can
be isolated in filtration, the only difference being that they are
lower (rather than higher) in filtration
grading than all other points.  In addition, two exterior points can be
removed by an isotopy.  The proof of Proposition
 \ref{prop:largen2} is exactly the same in this setting, with the
roles of subcomplex and quotient complex reversed (i.e.\ the cabled
knot group is naturally a quotient complex rather than a subcomplex,
and the first term in the short exact sequence of chain complexes has trivial homology). We state
the analogue of Theorem \ref{thm:largen} with $n<0$ for
completeness.   

\begin{theorem}
\label{thm:largenneg}
Let $K \subset S^3$ be a knot with $\mathrm{deg}~\HFKa(K)=d$. Then $\exists \ N<0$
and $c(c',n,p)$ such that $\forall \ n<N$, the following holds:  
$$ \mathrm{deg}~\HFKa(K_{p,pn+1})=pd+\frac{(p-1)(p|n|-2)}{2}. $$
If $i < -c(c',n,p)$ we have
{\small$$ \HFKa_*(K_{p,pn+1},i)\cong 
\left\{\begin{array}{ll}

 H_{*+2(d-k)}(\frac{\CFKa(K)}{\Filt(K,d-k-1)}) & {\text{$i=pk-pd-\frac{(p-1)(p|n|-2)}{2}$}} \\

 H_{*+2(d-k)-1}(\frac{\CFKa(K)}{\Filt(K,d-k-1)}) & {\text{$i=pk+1-pd-\frac{(p-1)(p|n|-2)}{2}$}} \\

 0 & {\text{otherwise.}}\\
\end{array}
\right. $$}%
When $p=2$ we can arrange that $c(c',n,p)<0$.
\end{theorem}

Proposition 3.7 of \cite{Knots} relates the Floer homology of a knot
$K$, and its mirror, ${\overline K}$.  Thus the above theorem,
together with Theorem \ref{thm:largen}, gives the values of
$\HFKa(K_{p,pn\pm 1},i)$ for all sufficiently large $|n|$ with $|i|>c(c',n,p)$.

\section{Examples}
\label{sec:examples}

Let $T_{2,2m+1}$ be the
$(2,2m+1)$ torus knot.  For all $i\ge0$ We have the following:

\begin{prop}
\label{prop:examples}

If $m>0$,$0\le k<m$, then for all $n>10m$
{\small
$$\HFKa((T_{2,2m+1})_{2,2n+1},i)\cong
\left\{\begin{array}{ll}
\Z_{-2k} & {\text{for $i=2m+n-4k$}} \\
\Z_{-2k-1} & {\text{for $i=2m+n-4k-1$}} \\
0 & {\text{for $i=2m+n-4k-2$}} \\
0 & {\text{for $i=2m+n-4k-3$}} \\
\Z_{i-n} & {\text{for $0\le i\le n-2m$}}\\
0 & {\text{otherwise}}\\
\end{array}
\right.$$}%
If $m<0$,$0\le k<|m|$, then for all $n>6|m|+1$
{\small
$$\HFKa((T_{2,2m-1})_{2,2n+1},i)\cong
\left\{\begin{array}{ll}
\Z_{2|m|-4k} & {\text{$i=2|m|+n-4k$}} \\
\Z_{2|m|-4k-1} & {\text{$i=2|m|+n-4k-1$}} \\
\Z_{2|m|-2k-1}\oplus \Z_{2|m|-4k-2} & {\text{$i=2|m|+n-4k-2$}} \\
\Z_{2|m|-2k-2}\oplus \Z_{2|m|-4k-3} & {\text{$i=2|m|+n-4k-3$}} \\
\Z_{i-n} & {\text{for $0\le i\le n-2|m|$}}\\
0 & {\text{otherwise}}\\
\end{array} 
\right.$$}%
If $m>0$,$0\le k<m$, then for all $n<-6m-1$
{\small
$$\HFKa((T_{2,2m+1})_{2,2n-1},i)\cong
\left\{\begin{array}{ll}
\Z_{-2m+4k} & {\text{$i=-2m-|n|+4k$}} \\
\Z_{-2m+4k+1} & {\text{$i=-2m-|n|+4k+1$}} \\
\Z_{-2m+2k+1}\oplus \Z_{-2m+4k+2} & {\text{$i=-2m-|n|+4k+2$}} \\
\Z_{-2m+2k+2}\oplus \Z_{-2m+4k+3} & {\text{$i=-2m-|n|+4k+3$}} \\
\Z_{i+|n|} & {\text{for $0\ge i\ge 2m-|n|$}}\\
0 & {\text{otherwise}}\\
\end{array}
\right.$$}%
If $m<0$,$0\le k<|m|$, then for all $n<-10|m|$
{\small
$$\HFKa((T_{2,2m-1})_{2,2n-1},i)\cong
\left\{\begin{array}{ll}
\Z_{2k} & {\text{for $i=-2|m|-|n|+4k$}} \\
\Z_{2k+1} & {\text{for $i=-2|m|-|n|+4k+1$}} \\
0 & {\text{for $i=-2|m|-|n|+4k+2$}} \\
0 & {\text{for $i=-2|m|-|n|+4k+3$}} \\
\Z_{i+|n|} & {\text{for $0\ge i\ge 2|m|-|n|$}} \\
0 & {\text{otherwise}}\\
\end{array}
\right.$$}
\end{prop}
\begin{proof}
When $n>0$ this follows directly from Proposition ~\ref{prop:T2n}
and Theorem ~\ref{thm:2ncable}.  When $n<0$ we use Theorem
~\ref{thm:largenneg}.  Alternatively, we could use Proposition 3.7 of
~\cite{Knots}. (In Proposition ~\ref{prop:T2n} we only
remarked upon the existence of certain differentials in the chain
complex for $S^3$.  Indeed, it is not difficult to see that these are
the only differentials.  Use Equations ~\eqref{eq:gr} and
~\eqref{eq:filt}, for instance, together with the fact that domains of
holomorphic disks must be positive.) 

 The bounds for $n$ when $n>0,m>0$ or $n<0,m<0$ arise as follows: In the Heegaard diagram for the $(2,2m+1)$
torus knot coming from its 2-bridge presentation, the natural
longitude has framing $n=4m$ (by natural we mean the longitude which
intersects $\alpha_2$ geometrically the minimum number of times).
Applying Lemma ~\ref{lemma:hd} with this longitude yields one point,
$x_0 \in \tilde{\beta}\cap\alpha_2$.  It is a straightforward
homological computation to see that when $m>0$ (resp.\ $m<0$), $x_0$ carries one intersection
point, $x_0\times w_{2m+1}$,
with relative filtration grading higher (resp. lower) than all other
points.  Furthermore, the
lowest (resp. highest) relative filtration dimension of any
intersection point in the
diagram is $12|m|$ less (resp.\ more) than $x_0\times w_{2m+1}$.
Thus the total breadth of relative filtration agrees with the breadth
of the Alexander polynomial and is equal to $12|m|+1$. By the lemmas in
Section ~\ref{sec:largen}, increasing (resp.\ decreasing) the
framing of the longitude by $6|m|$ increases the total breadth in
filtration by $12|m|$.  Thus the total breadth in filtration is $24|m|+1$
whereas the breadth in filtration of the exterior points is
$12|m|+1$. Therefore we may apply Theorem ~\ref{thm:2ncable}
(resp.~\ref{thm:largenneg}).  The bound in either case is $4m+6m=10m$. 

We determine the bound when  $n<0,m>0$ or $n>0,m<0$ as follows:  begin
as above with the natural $4m$-framed longitude for the $(2,2m+1)$ torus
knot. If $m>0$ (resp.\ $m<0$), the point $x_0\times w_1$ carried by $x_0$ with lowest
(resp.\ highest) relative filtration is
$4m$ lower (resp.\ higher) in filtration than $x_0\times w_{2m+1}$. 
This follows from Lemma ~\ref{lemma:varytuple}. Since the total breadth in filtration is
$12|m|+1$, we see that $x_0\times w_1$ is $8|m|$ higher (resp.\ lower)
in filtration
than the lowest (resp.\ highest) filtration dimension.  Decrease (resp.\ increase) the framing of the
longitude by $1$.  This has the effect of making the total breadth of
filtration $12m$, rather than $12m+1$. Furthermore, $x_0\times w_1$ is
now $8|m|-1$ higher
(resp.\ lower) than the lowest (resp.\ highest) filtration dimension.  This is because a pair of
points $x_1,x_2 \in \tilde{\beta}\cap\alpha_2$ can be removed by an
allowed isotopy as in Proposition ~\ref{prop:T2n}. The net effect
of the change in framing is moving the basepoint $w$ across
$\tilde{\beta}$.  Decrease (resp.\ increase) the framing by $4|m|$.
This ensures that $x_0\times w_1$  is lower in relative filtration than all other
points, interior or exterior. The framing used is now $-1$ (resp.\ $1$) and the total breadth of
filtration is $12|m|+1$.  Decrease (resp.\ increase) the framing
by $6|m|+1$.   The exterior points isolated in filtration
from the interior now account for more than half of the total filtration
breadth.  The bounds follow.
 \end{proof}

{\bf Remark}\qua It is likely that the bound here, and indeed in Theorem \ref{thm:2ncable}, can be taken to be $N=2g(K)$. Some support of this can be found in Chapter $3$ of \cite{mythesis}, where we explicitly compute all $(2,2n+1)$ cables of the trefoil knot.

\Addresses\recd
\end{document}